\documentclass{article} 

\usepackage{amssymb}

\begin{document}

\def \n {\noindent}
\def \R {{\mathbb R}} \def\C {{\mathbb C}}
\def \Q {{\mathbb Q}}
\def \Z{{\mathbb Z}}
\def \N {{\mathbb N}}
\def \ds {\displaystyle}
\def \bul {$\bullet$\ }
\def \ee {\prec} 
\def\cR {{\cal R}}
\def\cM{{\cal M}}
\def\cN{{\cal N}}
\def \L {{\cal L}}
\def \F {{\cal F}}
\def \exp{{\rm exp}}
%******************************************
\newtheorem{Theorem}{Theorem}[section]
\newtheorem{Proposition}[Theorem]{Proposition} 
\newtheorem{Lemma}[Theorem]{Lemma}
\newtheorem{Example}[Theorem]{Example}
\newtheorem{Corollary}[Theorem]{Corollary}
\newtheorem{Fact}[Theorem]{Fact}
\newtheorem{Conjecture}[Theorem]{Conjecture} 

\newenvironment{Definition} {\refstepcounter{Theorem} \medskip\noindent
 {\bf Definition
\arabic{section}.\arabic{Theorem}}\ }{\hfill}

\newenvironment{Remarks} {\refstepcounter{Theorem} 
\medskip\noindent {\bf Remarks
\arabic{section}.\arabic{Theorem}}\ }{\hfill}

\newenvironment{Exercise} {\medskip\refstepcounter{Theorem} \noindent {\bf Exercise
\arabic{section}.\arabic{Theorem}}\ } {\hfill}

\newenvironment{Proof}{{\noindent \bf Proof\ }}{\hfill}

\newenvironment{claim} {{\smallskip\noindent \bf Claim\ }}{\hfill}
%*******************************************

\title{Decidability of the Natural Numbers with the Almost-All Quantifier}

\author{David Marker\thanks{Partially supported by NSF grant DMS-0200393}\\
University of Illinois at Chicago \and
Theodore A. Slaman\thanks{Partially supported by NSF grant
  DMS-0501167.}\\University of California, Berkeley}

\maketitle

\begin{abstract}
  We consider the fragment $\F$ of first order arithmetic in
  which quantification is restricted to ``for all but finitely many.''
  We show that the integers form an $\F$-elementary substructure of the
  real numbers.  Consequently, the $\F$-theory of arithmetic is
  decidable.
\end{abstract}

\section{Introduction}

In this note, we consider the fragment of first order arithmetic in
which quantification is restricted to ``for all but finitely many''
and its negation ``there exist at most finitely many.''  These
quantifiers are quite natural, since mathematics is rich with deep
theorems asserting that the various potentially infinite sets are, in
fact, finite.  However, as we describe below, the formal quantifier
``for all but finitely many $x$'' is surprisingly weak.  In fact, the
fragment of arithmetic that it generates is unable to distinguish
between the natural numbers and the non-negative real numbers.
Consequently, we can use the quantifier elimination and decidability
of the theory of real closed fields, to deduce the same facts for this
fragment of the theory of arithmetic.

The proof that we give below is a direct application of the machinery
of o-minimality.  Even so, the question was originally motivated by
recursion theoretic investigations.  Typically, first order structures
are presented either by specifying their atomic diagrams or by
specifying their generators and relations.  The former is a recursive
($\Delta^0_1$) presentation of the structure and the latter is a
recursively enumerable ($\Sigma^0_1$) presentation.  The next level in
the arithmetic hierarchy is a $\Pi^0_1$ presentation of a first order
structure.  In a $\Pi^0_1$ presentation, the atomic truths are those
which are never canceled.  In particular, if one is attempting to give
a nonrecursive $\Pi^0_1$ presentation of the natural numbers, one is
faced with the problem of canceling an element $a$'s role as a
particular number and reassigning $a$ to an arbitrarily large value.
Asking about the atomic types of arbitrarily large numbers naturally
leads to asking about the almost-all theory of arithmetic.

No matter how the question happened to be asked, the model theory of
o-minimal structures provides a direct route to the answer.

\section{The Almost-All Theory of Arithmetic}

In the natural numbers the quantifier ``for all but finitely many'' is
equivalent to the quantifier ``for all sufficiently large''. Since we
will also consider non-discrete orderings it will be useful to work
with the later quantifier.

Let $\L=\{+,\cdot,<,0,1\}$ be the language of ordered rings.  Let $Q$
be a new quantifier symbol. If $\cM=(M,<,\dots)$ is a linearly ordered
structure, we say that
$$
 \cM\models Qx\ \phi\hbox{ if and only if } \cM\models \exists
 z\forall x>z\ \phi.
 $$
 Let $\F$ be the smallest collection of $\L$-formulas containing
 all quantifier free formulas and closed under propositional
 connectives and $Q$. We view $\F$ as a fragment of first order logic.

\begin{Theorem}\label{el-sub}
  The natural numbers is an $\F$-elementary substructure of the real
  ordered field.
\end{Theorem}

\begin{Proof} We will prove, by induction
  on complexity of $\F$-formulas, that if $\phi(x_1,\dots,x_m)$ is an
  $\F$-formula with free variables $x_1,\dots,x_m$ and
  $n_1,\dots,n_m\in \N$, then
$$\N\models \phi(n_1,\dots,n_m) \hbox{  if and only if  }
 \R\models \phi(n_1,\dots,n_m).$$
This is clear for atomic formulas and the induction is trivial for Boolean
combinations.

Suppose the claim is true for $\phi(x,\bar y)$.
If $\R\models Qx\ \phi(x,\bar n)$, then there is $r\in \R$
such that  $\R\models \phi(x,\bar n)$    for all $x>r$. If $s\in\N$
and $s\ge r$, then $\N\models\forall x>s\ \phi(x,\bar n)$.

On the other hand if $\N\models Qx\ \phi(x,\bar n)$,
then there is $r\in \N$ such that for all $s\in \N$ if $s>r$,
then $ \N\models \phi(s,\bar n)$ and $ \R \models \phi(s,\bar n)$.
Thus $X=\{x\in\R: \R\models\phi(x,\bar n)\}$ is unbounded.
But $\R$ is o-minimal and, hence, $X$ is a finite union of points
and intervals. Thus there is $r^\prime\in\R$ such that $(r,+\infty)\subseteq
X$ and $\R\models Qx\ \phi(x,\bar n)$.
\end{Proof}

\medskip In particular $\R$ and $\N$ have the same $\F$-theory.  We
can use the quantifier elimination and decidability of the theory of
real closed fields, to deduce the same facts for the $\F$-theory of
$\N$.

\begin{Corollary} i) The $\F$-theory of the natural numbers is decidable.

ii) Every $\F$-formula is equivalent in $\N$ to a quantifier free formula.
\end{Corollary}

 \medskip

The proof of Theorem \ref{el-sub} can be applied in more general settings.
Suppose $\cR=(\R,+,\cdot,<,\dots)$ is an 
o-minimal expansion of the real field in a language $\L$
and that the natural numbers are an $\L$-substructure of $\cR$.
Then $\N$ is an $\F$-elementary submodel of $\cR$. For example,
let $\cR=(\R,+,\cdot,<,e,0,1)$ where $e$ is the binary function
$$e(x,y)=\cases{x^y& if $x>0$\cr 0&otherwise\cr}.$$ The structure $\cR$
is a reduct of $\R_\exp$ the real field with exponentiation and 
Wilkie \cite{w} proved that $\R_\exp$ is o-minimal. Thus the natural numbers
is an $\F$-elementary submodel.
Macintyre and Wilkie \cite{mw} proved that if Schanuel's Conjecture is
true, then the theory of $\R_\exp$ is decidable.\footnote{Schanuel's Conjecture
asserts that if $\lambda_1,\dots,\lambda_n\in\C$ are $\Q$-linearly
independent, then the field $\Q(\lambda_1,\dots,\lambda_n,e^{\lambda_1},
\dots,e^{\lambda_n})$ has transcendence degree at least $n$ over $\Q$.}
Thus if Schanuel's Conjecture holds, 
then the $\F$-theory of $(\N,+,\cdot,<,x^y,0,1)$ is decidable.

\section{Completeness of the $\F$-Theory of Commutative Ordered Rings.}

Let $\L=\{+,\cdot,<,0,1\}$ be the language of ordered rings.
We conclude by proving that any two commutative ordered rings are
$\F$-elementarily equivalent.

 We always assume that rings have a multiplicative identity.

\begin{Definition} A commutative ordered ring $A$ is a {\em 
real closed ring} if the intermediate value property holds for
every polynomial in $A[X]$.
\end{Definition}

\medskip
Cherlin and Dickman \cite{cd} proved that real closed rings are exactly
convex subrings of real closed fields. 
For example, suppose $R$ is a real closed field containing
infinite elements. Then $A=\{x\in R: |x|<n$ for some $n\in\N\}$ is a real
closed ring that is not a field.
Cherlin and Dickman \cite{cd} also proved that any two real closed rings
that are not fields are elementarily equivalent.

Real closed rings need not be o-minimal.
In $(A,+,\cdot,<)$ we can define the monad of infinitesimals
as the ideal of noninvertible elements. This is a convex set that is
not an interval with endpoints in $A$.

\begin{Definition} A linearly ordered structure $(M,<,\dots)$ is 
{\em weakly o-minimal} if every subset of $M$ that is definable with
parameters from $M$ is a finite union of convex sets. 
\end{Definition}
\medskip

Dickman \cite{d} noted that all real closed rings are weakly o-minimal.
Indeed \cite{mms} shows that every weakly o-minimal ring is real closed.
Weakly o-minimal structures still satisfy the {\em Asymptotic Dichotomy
Principle}, namely, if $(M,<,\dots)$ is weakly o-minimal and $X\subseteq
M$ is definable, then there is $r\in M$ such that $(r,+\infty)\subseteq X$
or $(r,+\infty)\cap X=\emptyset$. This is enough to adapt the proof of 
Theorem \ref{el-sub} to prove a mild generalization.

\begin{Lemma}\label{gen}
Suppose  $A$ is a commutative ordered ring, $R$ is a
real closed ring and $A$ is unbounded in $R$.
If $\phi(x_1,\dots,x_n)$ 
is an $\F$-formula  with
free variables $x_1,\dots,x_n$ and $a_1,\dots,a_m\in A$, then
$$A\models \phi(a_1,\dots,a_m) \hbox{  if and only if  }
 R\models \phi(a_1,\dots,a_m).$$
\end{Lemma}

The assumption that $A$ is unbounded in $R$ is essential.
Suppose $F$ is a real closed field and $R$ is a proper convex
subring. Let $a\in R$ such that $1/a>R$. Then $R\models Qx\ ax<1$,
while this fails in $F$. While $R$ is not an $\F$-elementary submodel
of $F$, they will be $\F$-elementarily equivalent.

\begin{Theorem} Any two commutative ordered rings are $\F$-elementarily 
equivalent.
\end{Theorem}
\begin{Proof} Let $A$ be a commutative ordered ring. We will prove that
$A$ is $\F$-elementarily equivalent to $\Z$.

Let $F$ be a real closed field such that $A\subset F$ is bounded in $F$.
Let $R$ be the convex hull of $A$ in $F$ and let $R_1$ be the convex
hull of $\Z$. Since $R$ and $R_1$ are real closed rings that are
not fields they are elementarily equivalent.

Let $\phi$ be any $\F$-sentence, then
\begin{eqnarray*}
\Z\models\phi&\Leftrightarrow& R_1\models\phi, \hbox{ by Lemma \ref{gen}}\\
&\Leftrightarrow& R\models\phi\\
&\Leftrightarrow& A\models\phi \hbox{ by Lemma \ref{gen}}.
\end{eqnarray*}
\end{Proof}


\begin{thebibliography}{99}

\bibitem{d} M. A. Dickmann, Elimination of quantifiers for ordered valuation
rings,  J. Symbolic Logic (52) 1987, 116--128.

\bibitem{cd}G. Cherlin and M. A. Dickmann, Real closed rings II: model
theory, Ann. Pure Appl. Logic (25) 1983, 213--231.

\bibitem{mw} A. Macintyre and A. Wilkie, On the decidability of the real
 exponential field, {\em Kreiseliana: About and Around Georg Kreisel},
A K Peters, 1996, 441-467. 

\bibitem{mms}D. Macpherson, D. Marker, and C. Steinhorn, Weakly o-minimal 
structures and real closed fields, Trans. Amer. Math. Soc. (352) 5435--5483.

\bibitem{w} A. Wilkie, Model completeness results for expansions of the ordered field of real numbers by restricted Pfaffian functions and the exponential 
function, J. Amer. Math. Soc., 9(4) 1996, 1051-1094. 

\end{thebibliography}
\end{document}